\documentclass[MSNbibl,nameyear,seceqn,dvips]{arxstspdf}
\usepackage{flushend}
\usepackage{stfloats}
\volume{27}
\issue{4}
\pubyear{2012}
\firstpage{469}
\lastpage{480}
\doi{10.1214/12-STS397} 

\makeatletter
\newtheorem{theorem}{Theorem}[section]
\newproclaim{remark}{Remark}[section]
\newproclaim{example}{Example}[section]
\newproclaim{definition}{Definition}[section]
\newproclaim{condition}{Condition}
\newproclaim{KKTConditions}{KKT Conditions}

\newcommand{\R}{\mathbb{R}}
\newcommand{\PP} {\mathbb{P}}
\newcommand{\EE} {\mathbb{E}}
\makeatother

\begin{document}
\begin{frontmatter}

\title{Quasi-Likelihood and/or Robust Estimation in
High Dimensions}
\runtitle{Estimation in High Dimensions}

\begin{aug}
\author[a]{\fnms{Sara} \snm{van de Geer}\corref{}\ead[label=e1]{geer@stat.math.ethz.ch}}
\and
\author[b]{\fnms{Patric} \snm{M\"uller}\ead[label=e2]{mueller@stat.math.ethz.ch}}
\runauthor{S. van de Geer and P. M\"uller}

\affiliation{ETH Z\"urich}

\address[a]{Sara van de Geer is Full Professor, Seminar for Statistics, ETH Z\"
urich, R\"amistrasse 101, 8092 Z\"urich, Switzerland \printead{e1}.}
\address[b]{Patric M\"uller is Graduate Student, Seminar for Statistics, ETH Z\"
urich, R\"amistrasse 101, 8092 Z\"urich, Switzerland \printead{e2}.}

\end{aug}

%
\begin{abstract}
We consider the theory for the high-dimensional
generalized linear model with the Lasso.
After a short review on theoretical results in literature,
we present an extension of the oracle results to the case of
quasi-likelihood loss.
We prove bounds for the prediction error and $\ell_1$-error. The
results are derived under fourth moment conditions on the error distribution.
The case of robust loss is also given.
We moreover show that under an irrepresentable condition, the
$\ell_1$-penalized quasi-likelihood estimator has no false positives.
\end{abstract}

%
\begin{keyword}
\kwd{High-dimensional model}
\kwd{quasi-likelihood estimation}
\kwd{robust estimation}
\kwd{sparsity}
\kwd{variable selection}.
\end{keyword}

\end{frontmatter}
%

\section{A Review of the Theory in Literature} \label{review.section}

Consider $n$ independent observations
$\{ (x_i^T, Y_i) \}_{i=1}^n $, where $Y_i \in\mathcal{Y} \subset\R$ is
a random response
variable, and $x_i$ is a fixed $p$-dimensional vector of co-variables,
$i=1, \ldots, n$.
In a high-dimensional model, the number of co-variables $p$ is much
larger than the
number of observations $n$. There has been much literature on the
linear model
for this situation. In that case, one assumes that
\[
Y_i = x_i^T \beta^0 + \varepsilon_i ,\quad  i=1 , \ldots, n ,
\]
where $\beta^0 \in\R^p$ is an unknown vector of coefficients,
and $\varepsilon_1, \ldots, \varepsilon_n$ are independent noise
variables. The Lasso
estimator
(\cite{tibs96}) is
\[
\hat \beta:= \mathop{\operatorname{arg\,min}}_{\beta\in\R^p }   \Biggl\{ \sum_{i=1}^n |
Y_i -
x_i^T \beta|^2 + \lambda\sum_{j=1}^p | \beta_j |   \Biggr\} .
\]
The parameter $\lambda>0$ is a regularization parameter, and
$\| \beta\|_1 := \sum_{j=1^p} | \beta_j | $ is the $\ell_1$-norm of
$\beta$.
For the case of orthogonal design, that is, the case where the
columns of the $n \times p$ design matrix
\[
\mathbf{X}:=
\pmatrix{ x_1^T \cr \vdots\cr x_n^T \cr
}
\]
are orthogonal, the Lasso estimator is the soft-\linebreak[4]thresholding estimator
(\cite{donoho95}). We study in this paper the extension of the theoretical
results for the Lasso estimator, to the case of generalized linear models.

The theory for the Lasso with least squares loss is well established.
We refer
to \citet{Bunea:06}, \citet{Bunea:07a}, \citet{buneaetal06},
\citet{vandeG07},
\citet{Lou08}, \citet{bickel2009sal}. See also \citet{BvdG2011} and the
references
therein. The main results concern oracle inequalities for the
prediction error
$\| \mathbf{X} ( \hat\beta- \beta^0 ) \|_2^2$ and variable selection properties
of the Lasso. Oracle results say that the prediction error of the Lasso
estimator is up to log-factors as good as that of an oracle that uses
the least
squares ``estimator'' with only the co-variables in the unknown
active set $S_0 := \{ j \dvtx  \beta_j^0 \not= 0 \}$. Variable selection results
roughly state that with large probability, the estimated active set
$\hat S : = \{ j \dvtx  \hat\beta_j \not= 0 \}$ is with large probability
equal to the true active set~$S_0$. Both results depend on appropriate
conditions:
for prediction one assumes restricted eigenvalue condition (\cite
{koltch09a}; \cite{koltch09b}; \cite{bickel2009sal}) or compatibility
conditions (\cite{vandeG07}), and for variable selection, one assumes
the neighborhood
stability (\cite{mebu06}) or equivalent irrepresentable condition
(\cite{ZY07}).
Clearly, variable selection is a harder problem than prediction, so
that one
expects conditions for the former to be stronger than those for the latter.
Indeed, \citet{van2009conditions} show that the irrepresentable condition
implies the compatibility condition.

Concerning work on oracle inequalities for general loss,
an earlier paper which uses $\ell_1$-regularization in this context
is \citet{Loubes:02}. Here, the case of orthogonal
design is considered (thus, it has $p \le n $). The technique of proof
is, however,
very much along the lines of the later proofs for nonorthogonal design
(with possibly $p > n$), as developed by \citet{vandeG07} and others.
Some remarks on the proof technique can be found in \citet{vdG:2001},
highlighting
that with an $\ell_1$-penalty one can derive oracle inequalities with rates
faster than $1 / \sqrt n $, despite the fact that
the penalty-term $\lambda\| \beta^0 \|_1 $ itself is generally of
larger order than
$1/ \sqrt n $. The case of quantile regression was studied in
\citet{vandegeer03}, again only for the case of orthonormal design.
In \citet{Tarigan:06}, hinge loss with $\ell_1$-penalty is studied.
Here the design is not assumed to be orthogonal, and is in fact
random. This paper does not use restricted eigenvalue
or compatibility conditions, but rather a weighted eigenvalue
condition. It shows that the $\ell_1$-penalty leads to estimators
which are both adaptive to the ``smoothness'' or ``complexity'' of the
underlying regression function, as well as to the ``margin behavior''
of the problem. The margin behavior expresses the amount of
curvature of the theoretical risk near its minimum.
The paper
\citet{Bunea:07b} considers the density estimation problem.
In \citet{vandeG07}, results are derived for generalized
linear models with $\ell_1$-penalty and $p $ possibly larger than $n$,
assuming the
compatibility condition. It covers the case of quadratic loss and of
general Lipschitz loss,
and it allows for random design.
Similar results are in \citet{vandeG08}, although there the
compatibility condition
is replaced by one somewhat in the spirit conditions in \citet{JN11}.
In \citet{BvdG2011},\vadjust{\goodbreak} one can find further details concerning sparsity
oracle inequalities
for high-dimensional generalized linear models.

There is a large body of
literature extending the oracle results for the linear model to matrix
versions. It is beyond the scope of this paper to review this work,
and we only point to the generalization to robust loss, as given in
\citet{Candes09}.

Within this volume, the paper \citet{Wai11} gives a general account of
oracle results for
high-dimensional M-estimators. After our Theo-\break rem~\ref{oracle-quasi-likelihood.theorem},~we~briefly
discuss its
relation with \citet{Wai11}.

Concerning variable selection, the fact that the
irrepresentable condition is rather strong has led to
considering modifications of the Lasso, such as
two step procedures, and the SCAD introduced by \citet{Fan1997};
see,
for example, \citet{Wu09}
for the case of quantile regression.

Our paper focuses only on the theoretical aspects.
There is much literature on applications of the Lasso
in generalized linear models; see \citet{Wu:09}, for example. The
computational aspects are well-studied: see \citet{Friedman:10}.
The paper \citet{lambert2011robust} contains, apart from theory,
software descriptions and a real data example for the case of
Huber loss.
In \citet{wang2007robust}, $\ell_1$-regularization
with least absolute deviations loss is studied and
compared numerically with the least squares Lasso.

We present new results for prediction and variable selection for the
case of quasi-likelihood estimation. The findings for prediction are
along the lines
as those in \citet{vandeG08}, but this time completed with the
compatibility condition.
The paper details and extends the findings in \citet{BvdG2011}.
We also show that a weighted form of the irrepresentable condition
implies consistent variable selection.

\section{Quasi-Likelihood and Robust Loss} \label{introduction.section}

We model the dependence of the distribution of $Y_i$ on $x_i$ via a
linear function
$f_{\beta^0} (x_i) := x_i^T \beta^0$, where $\beta^0$ is a vector of unknown
coefficients. The problem is to estimate $\beta^0$ or the linear
predictor vector
$f_{\beta^0} := \mathbf{X} \beta^0 $, where $\mathbf{X}^T := (x_1, \ldots,
x_n )$.
We study a high-dimensional situation, where the number of variables
$p$ can be much larger than the sample size~$n$.
(For technical reasons, we assume that $p$ is at least~$2$.)
The vector $\beta^0$ is assumed to be sparse; that is, its number of nonzero
coefficients is assumed to be small. See Section~\ref{sparsity.section} for more
details on sparsity.

We consider two models. The first one is
a generalized linear model, with a given inverse link function $G$,
that is,
\[
\EE(Y_i \vert x_i) := \mu_0 (x_i) = G( x_i^T \beta^0),\quad  i=1 ,
\ldots, n ,
\]
with $\beta^0 \in\R^p$ a vector of unknown coefficients.
The quasi-(log)likelihood function is
\[
Q(y, \mu) := \int_{y}^{\mu} \frac{y-u}{V(u)} \,du ,\quad  y, \mu\in\mathcal{
Y} ,
\]
where $V\dvtx \R\rightarrow(0, \infty)$ is a given variance function;
see also \citet{mccullagh1989generalized}. Together, quasi-likelihood
and link function define quasi-likeli\-hood loss, as follows:

\begin{definition}\label{quasi-likelihood.definition}
The quasi-likelihood loss function is
\[
\rho(y, z ) := - Q( y , G(z)) ,\quad  y \in\mathcal{Y} ,  z \in\R.
\]
\end{definition}

In our second model, the dependence of the distribution of $Y_i$ on $x_i$
may be described through quantiles or other aspects of the
distribution. In particular, one can define
this dependence via
a loss function $\{ \rho(y, z)\dvtx  y \in\mathcal{Y},  z \in\R\} $, and
\[
f_i^0 := \arg\min_{z \in\R} \EE ( \rho(Y_i , z)
\vert x_i ) .
\]
The generalized linear model assumes
that $f_i^0 = x_i^T \beta^0 $ for some $\beta^0 \in\R^p$.

The robust case is the one where, for all $y\in\mathcal{Y}$, the loss
function $\rho(y, z)$ is Lipschitz
in $z$, with Lipschitz constant not depending on $y$. Without loss of generality
one can then assume the Lipschitz constant to be equal to one.
This leads to the following definition:

\begin{definition}\label{robust.definition}
The loss function $\rho$
is \textit{robust} if for all $y \in\mathcal{Y}$,
\[
| \rho(y, z) - \rho(y, \tilde z)| \le|z - \tilde z | \quad \forall z,
\tilde z .
\]
\end{definition}

Quasi-likelihood loss is sometimes robust, but there are also many examples
where it is not. Moreover, there are many (robust) loss functions which
do not correspond to minus quasi-likelihoods.
See Section~\ref{examples.section} for some
examples.

To handle the large $p$ situation, one needs a regularized estimation method.
Let us write a linear function with coefficients $\beta$ as
\[
f_{\beta} (x) = x^T \beta.
\]
In what follows, we sometimes, with some abuse of notation, let
$f_{\beta}$ be the $n$-dimensional vector $\mathbf{X} \beta= (f_{\beta}
(x_1) , \ldots,
f_{\beta} (x_n))^T \in\R^n$ as well.

The $\ell_1$-norm of a vector $\beta\in\R^p$ is
\[
\| \beta\|_1 := \sum_{j=1}^p|\beta_j| .
\]
We examine the $\ell_1$-penalized estimator $\hat\beta$ of $\beta^0$,
defined as
\[
\hat\beta:=
\arg\min_{\beta\in\R^p}  \Biggl\{\frac{1}{n}
\sum_{i=1}^n \rho( Y_i, f_{\beta} (x_i) ) + \lambda\| \beta\|_1
 \Biggr\} .
\]
Here, $\lambda>0$ is a tuning parameter. Large values
correspond to more regularization, which means more shrinkage of the estimator
$\hat\beta$. The expression
\[
{1 \over n} \sum_{i=1}^n \rho( Y_i, f_{\beta} (x_i) )
\]
is called the empirical risk (at $\beta$). For least squares loss (i.e.,
$\rho(y, u) = (y-u)^2 $), the empirical risk is the usual
sum of squares (normalized by $1/n$). The above estimator is then
called the Lasso
estimator (\cite{tibs96}).

We will study loss functions $\rho$ that are either minus
quasi-likelihoods or robust
(or both).
The normalized Euclidean norm on $\R^n$ is
\[
\| f \|_n :=\sqrt{ f^T f / n } ,\quad  f \in\R^n .
\]
We will establish bounds for the
``prediction error'' $\| f_{\hat\beta} - f_{\beta^0}\|_n^2$, the
$\ell_1$-error
$\| \hat\beta- \beta^0\|_1$, and (for the case of
quasi-likelihood loss) present sufficient conditions for variable selection
using $\hat\beta$.

\subsection{Convex Loss} \label{convex.section}
We require throughout this paper, both for quasi-likelihood loss as well
as for robust loss, that the map
\[
z \mapsto\rho( y, z)
\]
is convex for all $y \in\mathcal{Y}$. This assumption is important from a
computational
point of view. It also plays a crucial role in our theory, as it allows us
to prove that the estimator $\hat\beta$ is in an $\ell
_1$-neighborhood of
$\beta^0$. This in turn will be invoked to establish sup-norm bounds
for $f_{\hat\beta}$.

\subsection{Sparsity}\label{sparsity.section}
The indices of the set of nonzero coefficients of $\beta^0$ is called
the (true) active set. It is denoted by
\[
S_0 := \{ j\dvtx  \beta_j^0 \not= 0 \} .
\]
Its cardinality $s_0 := | S_0 |$ is called the sparsity index of~$\beta^0$.
It is assumed that $s_0$ is relatively small, at least smaller than
$\sqrt{n/\log p}$ in order of magnitude; see (\ref{s0}), (\ref{s02}),
(\ref{s03}), (\ref{s04}) and (\ref{s05}). The vector $\beta^0$ is
sparse if
$s_0$ is small.

More generally, one can call a vector $\beta^0$ sparse if it can in some
sense be approximated by a vector with only a few nonzero entries.
To avoid too many digressions, we will not elaborate on this issue,
but only present a brief outline after the formulation of the main
oracle result; see Remark~\ref{oracle-quasi-likelihood.remark}.

\subsection{Results in this Paper} \label{results.section}
As $\beta^0$ is unknown, its active set $S_0$ and its sparsity index
$s_0$ are
unknown as well. We will show in Theorems~\ref{oracle-quasi-likelihood.theorem}
and~\ref{oracle-robust.theorem} that the prediction error of the
$\ell_1$-penalized estimator $\hat\beta$ is, up to a $\log p$-term,
the same
as that of minimizer of the empirical risk without penalty, but with
all coefficients not in
$S_0$ restricted to be zero. The latter is not an estimator, as it depends
on the unknown $S_0$. It is often referred to as the oracle.
We moreover show that a version of the irrepresentable condition,
appropriate for
quasi-likelihood loss, is sufficient for variable selection; see
Theorem~\ref{selection-quasi-likelihood.theorem}.
All our results are stated in a nonasymptotic form, but to facilitate
the interpretation, we also
give asymptotic formulations.

%
%
%

\subsection{Organization of the Paper} \label{organization.section}
The next section provides some examples of quasi-likelihood and robust loss.
Section~\ref{compatibility.section} gives the definition of the so-called
\textit{compatibility constant}, which will occur in the oracle results.
Section~\ref{oraclequasi-likelihood.section} gives oracle inequalities
for the prediction and $\ell_1$-error for quasi-likelihood loss, and Section
\ref{oracle-robust.section} does the same for robust loss.
In Section~\ref{selection-quasi-likelihood.section} we address the variable
selection problem in the quasi-likelihood context. Similar arguments
can be used in the robust context, but this is omitted here.
Section~\ref{random-design.section} briefly discusses the case of
random design,
and Section~\ref{discussion} concludes.
The proofs are in the supplemental article \citet{GeerMullersupp:12}.
Lemmas A.2 and A.7 there are
based on a concentration inequality (see \cite{Massart:00a}) and a
contraction inequality; see
\citet{Ledoux:91}. These lemmas use only fourth moment assumptions, and
are perhaps of
interest in themselves.

\section{Examples of Loss Functions}\label{examples.section}

\subsection{Least Squares Loss}
The least squares criterion has $\mathcal{Y} = \R$. It corresponds to a
quasi-likelihood loss with variance function
$V(u) =1$ for all $u \in\R$. The link function is then the identity,
which is
the canonical link function for this case.
The loss function is convex, but not robust.

\subsection{Logistic Loss}\label{logistic.section}
When the response $Y_i$ is binary, say $Y_i \in\{ 0, 1 \} $, $i=1 ,
\ldots, n $,
we have
\[
\EE(Y_i \vert x_i) = \PP( Y_i =1 \vert x_i) .
\]
In logistic regression, one takes the quasi-likelihood with variance function
$V(u) = u(1-u)$, $u \in(0,1)$, and the canonical link function
\[
\gamma(\mu) := \log \biggl( \frac{\mu}{1- \mu}  \biggr) ,\quad  \mu\in
(0,1),
\]
that is,
\[
G(z) = \gamma^{-1} (z) = \frac{\mathrm{e}^z}{1+ \mathrm{e}^z } ,\quad  z \in\R
.
\]
Hence, in this case,
\[
\rho(y,z) = y z - \log(1+ \mathrm{e}^z) ,\quad  z \in\R.
\]
Because $\mathcal{Y}= \{ 0,1\} $, one sees that this leads to a robust loss
function, that is, $ z \mapsto\rho(y,z) $ is Lipschitz in $z$ for all
$y \in\mathcal{Y}$.
We acknowledge that logistic regression is not robust in the sense
of having a bounded influence function (but we will in fact assume in
Condition~\ref{conA1} that the covariables are bounded).
As in all cases of quasi-likelihood with canonical link
function, the loss also convex.

\subsection{Binary Response with Other Link Functions}\label{binary.section}
Consider binary response $Y_i \in\{ 0, 1 \}$ as in Section~\ref{logistic.section}, but now with more general inverse link
function $G$.
\[
\PP(Y_i = 1 \vert x_i) = G (x_i^T \beta^0 ) ,\quad  i=1 , \ldots, n .
\]
If $G\dvtx \R\rightarrow[0,1] $ is a strictly increasing symmetric
distribution function, then quasi-likelihood loss is convex. This is because
the hazard $g(u)/ (1- G(u))$ ($g$ being the derivate of $G$) is a
decreasing function of~$u$. When the hazard
is uniformly bounded, quasi-likelihood loss is also robust.

\subsection{Quantile Regression}\label{quantile.section}
If the dependence of the distribution of $Y_i \in\R$ on $x_i$ is via its
$\alpha$-quantile ($0 < \alpha< 1$), we take as loss function
\[
\rho(y, z) = \rho(y-z),
\]
where
\[
\rho(z) = \alpha| z| \mathrm{l} \{ z > 0 \} + (1- \alpha) |z| \mathrm{l}
\{ z \le0 \} .
\]
This is clearly a robust loss function, but it does not correspond
to a quasi-likelihood.

\section{The Compatibility Condition} \label{compatibility.section}

Let $S\subset\{ 1 , \ldots, p\}$ be an index set with cardinality~$s$.
We define for all $\beta\in\R^p$,
\[
\beta_{S,j} := \beta_j \mathrm{l} \{ j \in S\} ,\quad  j=1 , \ldots, p,\quad
\beta_{S^c} := \beta- \beta_S .
\]

Below, we present for constants $L>0$ the compatibility constant $\phi(L,S)$
introduced in \citet{vandeG07}.
For normalized design (i.e., $\| \mathbf{X}_j \|_n = 1 $ for all~$j$,
where $\mathbf{X}_j$
denotes the $j$th column of $\mathbf{X}$), one
can view $1- \phi^2 (1,S)/2$ as an $\ell_1$-version of the canonical
correlation
between the linear space spanned by the variables in $S$ on the one hand,
and the linear space of the variables in $S^c$ on the other hand.
Instead of all linear combinations with normalized $\ell_2$-norm,
we now consider all linear combinations with normalized $\ell_1$-norm
of the coefficients.
For a geometric interpretation, we refer to \citet{Lederer:11}.

\begin{definition*}
 The \textit{compatibility constant} is
\[
\phi^2 (L, S) := \min\{s \| f_{\beta} \|_n^2 \dvtx  \|\beta_S \|_1 =1
,
\| \beta_{S^c} \|_1 \le L \} .
\]
\end{definition*}

The compatibility constant is closely related to (and never smaller
than) the restricted
eigenvalue as defined in \citet{bickel2009sal}, which is
\[
\phi_\mathrm{RE}^2 (L,S)=
\min \biggl\{ { \| f_{\beta} \|_n^2 \over
\| \beta_S \|_2^2 } \dvtx \|\beta_{S^c} \|_1 \le L
\| \beta_{S} \|_1   \biggr\} .
\]

The calculation of the compatibility constant is a nonlinear
eigenvalue problem
[see, e.g., \citet{Hein:2010} for computational aspects of nonlinear
eigenvalues].
Lower bounds that hold with high probability follow, for example, if
$\mathbf{X}$ is an i.i.d. sample from
a $p$-dimensional vector with nongenerate covariance matrix; see
Section~\ref{random-design.section}
for some details.
See also \citet{koltch09a}, and see \citet{van2009conditions} for a
discussion of the
relation between restricted eigenvalues and compatibility.

For oracle results, we need $\phi(L, S_0)$ to be strictly positive
for some $L>1$ (depending on the tuning parameter $\lambda$).
In this paper, we take $L=3$ for definiteness, and we require
throughout that $\phi(3, S_0)>0$ (except when
we consider sparse approximations of the truth; see Remark
\ref{oracle-quasi-likelihood.remark}). If $\phi(3, S_0) =0$, one sees
that some conditions
[e.g., condition (\ref{s0})] become
impossible.

As we will see, all bounds in this paper involve not so much the
sparsity index $s_0$ itself,
but rather the \textit{effective sparsity}
\[
\Gamma_\mathrm{effective} (S_0) := {s_0 \over\phi^2 (3, S_0) } .
\]

\begin{example}
As a simple numerical example, let us suppose $n=2$,
$p=3$, $S_0= \{ 3 \}$ and
\[
\mathbf{X} = \sqrt{n}
\pmatrix{ 5/13 & 0 & 1 \cr12/13 & 1 & 0 \cr
}
.
\]
Thus, the sparsity index is $s_0=1$.
One can easily verify that there is no $\beta\in\R^p$ with
$\mathbf{X} \beta= 0 $ and $\| \beta_{S_0^c} \|_1 \le3 \| \beta_{S_0}
\|_1 $. Thus,
the compatibility constant $\phi^2(3, S_0) $ is strictly positive. In
fact, $\phi(3,S_0)$ is
equal to the distance of $\mathbf{X}_1$ to line that connects $3 \mathbf{X}_1
$ and $-3 \mathbf{X}_2 $,
that is $\phi(3,S_0) = \sqrt{2/13} $. The effective sparsity is
$\Gamma_\mathrm{effective} (S_0)=13/2$.

Alternatively, when
\[
\mathbf{X} = \sqrt{n}
\pmatrix{ 12/13 & 0 & 1 \cr5/13 & 1 & 0 \cr
}
,
\]
then $\phi(3,S)=0$. This is due to the sharper angle between $\mathbf{
X}_1$ and $\mathbf{X}_3$.
\end{example}

\section{Oracle Inequalities for Quasi-Likelihood Loss}\label
{oraclequasi-likelihood.section}

\subsection{The Case of Least Squares Loss}\label{leastsquares.section}
To appreciate the results we will present for the general case, it may
be useful
to first reconsider the standard linear model and least squares loss.
Let $Y = ( Y_1 , \ldots, Y_n)^T $ and suppose
\[
Y = \mathbf{X} \beta^0 + \varepsilon.
\]
Let $\hat\beta$ be the Lasso estimator
\[
\hat\beta= \arg\min_{\beta\in\R^p }  \{ \| Y- \mathbf{X}
\beta\|_n^2 + \lambda\| \beta\|_1  \} .
\]
Let $\mathbf{X} _j$ denote the $j$th column of the design matrix~$\mathbf{
X}$. If the errors $\varepsilon;=
(\varepsilon_1 , \ldots, \varepsilon_n )^T $ are independent with mean
zero and the design is normalized (i.e., \mbox{$\| \mathbf{X }_j \|_n =1$} for
all $j$) one can prove that uniformly in $j$, the ``correlations''
$\varepsilon^T X_j/n $ are small in absolute value,
generally as small as $O (\sqrt{\log p / n })$. The regularization
parameter $\lambda$ is to be chosen
in such a way that it ``overrules'' these correlations.
Indeed, this allows one to prove the following result
[see \citet{BvdG2011}, Theorem~6.1] by rather elementary means (recall
the notation
$f_{\beta} := \mathbf{X} \beta$):

\begin{theorem}\label{leastsquares}
Suppose that
\[
\lambda\ge4 \max_{1 \le j \le p } | \varepsilon^T  \mathbf{X} _j |/n .
\]
 Then
\[
\| f_{ \hat\beta} - f_{\beta^0} \|_n^2 + \lambda\| \hat\beta-
\beta^0 \|_1 \le
4 \lambda^2 \Gamma_\mathrm{effective } (S_0) .
\]
\end{theorem}

This result says that if the effective sparsity\linebreak[4] $\Gamma_{\mathrm{effective}} (S_0) $
is of the same order as the sparsity index $s_0 := |S_0|$ (i.e., if the
compatibility constant stays away from zero), then for a large class
of error distributions the Lasso estimator with
$\lambda\asymp\sqrt{\log p / n }$ is up to constants and a $(\log
p)$-factor as
good as as the oracle least squares ``estimator'' which knows the active
set $S_0$.
The performance of $\hat\beta$ is here measured
in terms of its prediction error\footnote{The prediction error
of the predictor $f_{\hat\beta} $ of an independent copy $Y_\mathrm{new}
:= f_{ \beta^0} + \varepsilon_\mathrm{new}$
of $Y$ is rather $\| f_{\hat\beta} - f_{\beta^0 } \|_n^2 + \sigma^2
$, where
$\sigma^2 = \EE\| \varepsilon_\mathrm{new} \|_n^2 $. We however do not
include the additional variance
$\sigma^2$ in our definition.} $\|X( \hat\beta- \beta^0 ) \|_n^2 $.
Theorem~\ref{leastsquares} moreover says that the $\ell_1$ error converges
with rate $\lambda\Gamma_\mathrm{effective} (S_0)$. Looking ahead at
more general
loss functions, ideas are based on quadratic approximations, which are
generally only valid in a neighborhood of $\beta^0$. This is why in
our work,
we will assume that $\lambda\Gamma_\mathrm{effective} (S_0)$ is small,
say $ \lambda\Gamma_\mathrm{effective} (S_0) \le\gamma$, where $\gamma
$ is
a sufficiently small constant.
With $\lambda\asymp\sqrt{\log p / n}$, and a compatibility constant
staying away from zero, it means we assume the sparsity index
$s_0$ to be sufficiently smaller than $\sqrt{n / \log p}$.

\subsection{General Quasi-Likelihood Loss}
As in the situation of the standard linear model and least squares loss,
we will study the error $\| f_{\hat\beta} - f_{\beta^0} \|_n^2 $ and the
$\ell_1$-error. For prediction, one will be interested in
estimating the mean $\mu_0 = G( f_{\beta^0})$ of the
response variable $Y$. Our Conditions~\ref{conA3} and~\ref{conA4} below will ensure that
$G$ has a bounded derivative on an appropriate domain.
This means that bounds for $\| f_{\hat\beta} - f_{\beta^0 }\|_n$
immediately lead to similar bounds for
$\|G( f_{\hat\beta} ) - G( f_{\beta^0 } ) \|_n$.
With some abuse of terminology, we refer to $\| f_{\hat\beta} -
f_{\beta^0 } \|_n^2$ as
the prediction error.\looseness=1

The theoretical properties of the $\ell_1$-penalized qua\-si-likelihood estimator
$\hat\beta$ depend on the tail-behavior of the error
\[
\varepsilon_i := Y_i - \mu_0 (x_i) ,\quad  i=1 , \ldots, n .
\]
We will need at least finite second moments of the errors.
For definiteness, we assume the errors have finite fourth moments.
With higher order moments, the confidence level
in the oracle result of Theorem~\ref{oracle-quasi-likelihood.theorem} will be larger, and when the errors
have sup-exponential tails, one can derive exponential probability
inequalities for prediction error and $\ell_1$-error.

\renewcommand{\thecondition}{$\mathrm{A}_{\varepsilon}$}
\begin{condition}\label{conA}
There exist constants
$\sigma>0$ and $\kappa>0$ such that
\[
\max_{1 \le i \le n } \EE\varepsilon_i^2 \le\sigma^2
\]
and
\[
{1 \over n} \sum_{i=1}^n \EE ( \varepsilon_i^2 - \EE\varepsilon
_i^2  )^2 \le
\kappa^4 .
\]
\end{condition}

The next conditions,~\ref{conA1}--\ref{conA4}, allow us to use quadra\-tic
approximations in a neighborhood of $\beta^0$.
We assume throughout that the inverse link function $G$ is increasing
and that its derivative
\[
g(z) := \frac{ d G(z)}{dz } ,\quad  z \in\R,
\]
exists.
We further define
\begin{eqnarray}\label{gammaB}
\gamma(\mu ) &:=& \int_{y_0}^{\mu } \frac{1}{V(u) }\, du ,\nonumber\\[-8pt]\\[-8pt]
B(\mu, \mu_0 ) &:=& \int_{\mu_0}^{\mu} \frac{ u- \mu_0}{V(u) }\, du ,\quad
 \mu\in\mathcal{Y} ,\nonumber
\end{eqnarray}
where $y_0$ is an arbitrary but fixed constant.
We let
%
\begin{equation}\label{H}
H(z) := \gamma(G(z)) ,\quad  z \in\R,
\end{equation}
that is, $H := \gamma\circ G$.
Note that $\gamma$ is (up to an additive constant) the canonical link
function. When $G=\gamma^{-1}$, we get $H(z) = z$ for all $z$.
The term $y H(z)$ in the quasi-likelihood $Q(y, G(z))$ containing the
response $y$ is then linear in $z$.
In a sense, $H$ measures the departure from linearity of this term.
We let
\[
h(z) :=\frac { d H(z)}{dz } = {g(z) \over V(G(z))} ,\quad  z \in\R.
\]
The quantity $B(\mu, \mu_0)$ is the ``regret'' for choosing the
expectation $\mu$
instead of the ``true'' $\mu_0$.

\renewcommand{\thecondition}{\textup{A}\arabic{condition}}
\setcounter{condition}{0}
\begin{condition} \label{conA1}
There exists a constant $K_X$ such that
\[
\max_{1 \le j \le p } \max_{1 \le i \le n } | x_{i,j} | \le K_X .
\]
\end{condition}

We remark that Condition~\ref{conA1} serves as normalization of
the design, albeit not in terms of the $\| \cdot\|_n$ norm
but rather in supremum norm.
As our results will be presented in nonasymptotic form,
it is in principle possible to see the effect when, say,
$K_X$ grows with $p$ and/or $n$.

\begin{condition} \label{conA2}
There exists a constant $K_0$ such that
\[
\max_{1 \le i \le n } | f_{\beta^0 } (x_i) | \le K_0 .
\]
\end{condition}

\begin{condition} \label{conA3}
With $K_X$ and $K_0$ given in Conditions~\ref{conA1} and~\ref{conA2}
respectively,
there exists a positive constant $C_h$ such that for all
$|z| \le K_x + K_0 $,
\[
1 / C_h \le{ h(z) } \le C_h .
\]
\end{condition}

\begin{condition} \label{conA4}
With $K_X$ and $K_0$ given in Conditions
\ref{conA1} and~\ref{conA2} respectively, there exists a constant~$C_V$,
such that for all $|z| \le K_X + K_0$,
\[
2/ C_V \le V\circ G (z) \le C_V/2 .
\]
\end{condition}

\begin{remark}\label{examples.remark}
There is an interplay between Conditions~\ref{conA},~\ref{conA1} and~\ref{conA2}.
For example, for quadratic loss, we do not need~\ref{conA1} and~\ref{conA2} when the errors
are (sub)Gaussian. Conditions~\ref{conA1} and~\ref{conA2} are imposed so that we
need Conditions~\ref{conA3} and~\ref{conA4} only in the neighborhood
$|z| \le K_X + K_0 $. As for Condition~\ref{conA3},
when $G$ is the inverse of the canonical link function $\gamma$, it
holds with $C_h=1$, as $H$ is then the identity. For quadratic loss,
and logistic loss for example (which have canonical link function),
Condition~\ref{conA4} holds as well. We actually will only need the
lower bound for $V \circ G$ in this section, and the upper bound will come
into play in Section~\ref{selection-quasi-likelihood.section}.
\end{remark}

To organize the constants appearing in our results, let use the short
hand notation
\begin{eqnarray*}
C_{h,V}&:=& C_V C_h^2 ,
\\
C_{h,X}&:=& 16 C_h K_X,
\\
\Gamma(S_0) &:=& 16 C_{h,V} \Gamma_\mathrm{effective} ( S_0) .
\end{eqnarray*}
Thus, up to constants $\Gamma(S_0)$ is the effective sparsity.
As in the case of least squares loss, we assume the regularization
parameter $\lambda$ to be of order at least $\sqrt{\log p / n }$.
The larger $\lambda$,
the larger the confidence level of our bounds will be (in Theorem
\ref{leastsquares} this the probability of $4 \max_{1 \le j \le p } |
\varepsilon^T \mathbf{X}_j |/ n \le\lambda$)
but then these bounds themselves are also larger.
We introduce a variable $t>0$ to describe this effect, and
define
\[
\lambda_{\varepsilon} (t) :=
C_{h,X} \sigma\sqrt{2(t + \log p) \over n} .
\]
If we choose the tuning parameter $\lambda$ at least as
large as $4 \lambda_{\varepsilon} (t)$, the
confidence level will be at least $1- \alpha(t)$, where
\[
\alpha(t) :=
\alpha(t):= 3\exp[-t] + 3 \kappa^4 /(n \sigma^4).
\]
The variable $t$ is in principle arbitrary,
but it is, however, not allowed to be arbitrarily large.
As we can only apply the quadratic approximations in
a neighborhood of $\beta^0$ we will need to show that
$\hat\beta$ is with large probability\vadjust{\goodbreak} in such a neighborhood.
For that reason, we cannot let the tuning parameter $\lambda$
to be arbitrarily large (as a large $\lambda$ will give slow rates);
see condition (\ref{lambda.equation}) in Theorem
\ref{oracle-quasi-likelihood.theorem} below.
A~reasonable choice for $t$ is for example $t \asymp\log n $,
in which case $\alpha(t) \asymp1/n $.

\begin{theorem} \label{oracle-quasi-likelihood.theorem}
Let $\hat
\beta$ be the
$\ell_1$-penalized quasi-likelihood estimator.
Assume Conditions
\ref{conA} and
\ref{conA1}--\ref{conA4}.
Suppose that
%
\begin{equation} \label{s0}
\lambda_{\varepsilon} (t) \Gamma(S_0) \le\tfrac{1}{4 } .
\end{equation}
Take
%
\begin{equation}\label{lambda.equation}
4 \lambda_{\varepsilon} (t) \le\lambda\le{1 \over\Gamma(S_0) } .
\end{equation}
With probability at least $1- \alpha(t)$, it holds that
\[
\| \hat\beta- \beta^0 \|_1 \le{\lambda\over2} \Gamma(S_0)
\]
and
\[
\| f_{\hat\beta} - f_{\beta^0} \|_n^2 \le
\tfrac{ 3}{4} C_{h,V} \lambda^2 \Gamma(S_0) .
\]
\end{theorem}

\begin{remark}\label{wainwright}
Our result in Theorem~\ref{oracle-quasi-likelihood.theorem} is
comparable to
Corollary 3 in \citet{Wai11}, albeit that we do not assume bounded
responses or canonical link function, and
our compatibility condition is weaker than their assumed restricted eigenvalue
condition. On the other hand, we require (\ref{s0}), and only give
bounds for the
$\ell_1$-error and prediction error, not for the
$\ell_2$-error.
\end{remark}

\begin{remark}\label{constants}
We have presented the result in a nonasymptotic form,
but did not try to optimize the constants.
\end{remark}

\begin{remark}\label{interpretation.remark}
Thus, up to the compatibility constant, and taking $\lambda$ of order
$\sqrt{\log p / n }$, the prediction error is
of order $s_0 \log p / n $.
\[
\| \hat f - f^0 \|_n^2 = {\mathcal O}  \biggl( {s_0 \log p \over n }
 \biggr) .
\]
An oracle that knows $S_0$ and
does empirical risk minimization without penalty but with the
restriction that all coefficients not in $S_0$ are set to zero,
has a prediction error of order $s_0/n$. We see that for not knowing~$S_0$, one pays a price of order $\log p $.
We moreover have
\[
\| \hat\beta- \beta^0 \|_1 = {\mathcal O}  \Biggl( s_0 \sqrt{\frac{\log p
}{ n }}  \Biggr) .
\]
\end{remark}

\begin{remark}\label{oracle-quasi-likelihood.remark}
We have presented the above oracle inequality involving the sparsity
of the true $\beta^0$. If the\vadjust{\goodbreak} truth is not sparse, or if actually the
generalized linear model is misspecified, one may replace the
truth by a sparse linear approximation of the truth, and the
oracle inequality involves a trade-off between the approximation error
on the one hand, and the sparsity and compatibility constant on the other.
This trade-off is of the following form. Let for an arbitrary index
set $S \subset\{ 1 , \ldots, p\}$,
\[
\mathrm{f}_S := \arg\min_{f = f_{\beta_S}} \bar B_n ( G \circ f , \mu
_0 ) ,
\]
where $\bar B ( G \circ f , \mu_0)$ is the average regret
\[
\bar B ( G \circ f , \mu_0) := \frac{1}{n}
\sum_{i=1}^n B\bigl(G\circ f (x_i) , \mu_0 (x_i)\bigr) .
\]
Thus, $\mathrm{f}_S $ is the best approximation of $f^0$ using only the
variables in $S$.
Then under some regularity conditions the prediction error of $\bar B(
G \circ f_{\hat\beta} ,
\mu_0)$ is with
probability $(1-\alpha)$ bounded by
\[
\mathrm{const.} \min_{\mathrm{sets}  S}  \biggl\{
\bar B(G\circ\mathrm{f}_S , \mu_0 ) + \frac{\lambda^2 |S| }{
\phi^2 (L, S)}  \biggr\} .
\]
The ``const.'' depends on the constants occurring in the regularity conditions,
the constant $L$ depends moreover on the choice of $\lambda$, and
the confidence level $\alpha$ depends on all these.
For more details on this extension, we refer to \citet{BvdG2011} and
the references therein.
\end{remark}

\begin{remark} \label{s0.remark}
Condition (\ref{s0}) assumes that the sparsity index $s_0$ is
sufficiently smaller
than\break  $\sqrt{n / \log p }$, a condition we already announced in Section
\ref{leastsquares.section}.
This assumption plays its part in all our results: it
will also be important for variable selection and simplifies
the derivation of results for the case of random design.
In the case of least squares loss, the assumption can be avoided,
even in some cases with random design. It should, however, be noted that
a large $s_0$ means a slow rate. In particular, when
the sparsity is of larger order than $\sqrt{n / \log p }$, the bound for
the prediction error is of larger order than $\sqrt{\log p / n }$, and this
cannot be improved up to the $\log p$-term. Thus, then the bounds
are actually quite large in order of magnitude. Indeed, recall that the
prediction error is
$\| f_{\hat\beta} - f_{\beta^0} \|_n^2$, which is the squared
distance between $f_{\hat\beta}$ and $f_{\beta^0 } $.
Assumption~(\ref{s0}) allows us to conclude that
$\| \hat\beta- \beta^0 \|_1 \le1$, and hence, that $|f_{\hat\beta
} (x_i) | \le K_X + K_0 $
for all $i$. The latter was used because we only want to require
Conditions~\ref{conA3} and~\ref{conA4}
for bounded values of the argument $z$. When dealing with least squares loss,
Conditions~\ref{conA3} and~\ref{conA4} hold for all $z \in\R$. This means that with
least squares loss,
Assumption (\ref{s0}) can be dropped in Theorem \ref
{oracle-quasi-likelihood.theorem}; see Theorem~\ref{leastsquares}.
\end{remark}

\begin{remark}\label{sigma.remark}
The lower bound in (\ref{lambda.equation}) for the tuning parameter
$\lambda$ depends on the noise level
$\sigma$ as well as other unknown constants. In practice, one may for
instance apply cross-validation.
The noise level $\sigma$ can also be treated as additional parameter
which can be estimated
along with $\beta^0$. See \citet{Staedler:10} for a discussion.
\end{remark}

\section{Oracle Inequalities for Robust Loss}\label{oracle-robust.section}
In this section, we assume throughout that $\rho$ is robust loss;
see Definition~\ref{robust.definition}.

We define for $i=1, \ldots, n $,
\[
l_i (z) = \EE\rho(Y_i , z \vert x_i) ,\quad  z \in\R
\]
and assume
that $\ddot l_i (z) := d^2 l_i (z)/ dz^2 $ exists.

\renewcommand{\thecondition}{\textup{B}}
\begin{condition} \label{conB}
For $K_X$ and $K_0$ given in Conditions~\ref{conA1} and~\ref{conA2}
respectively, we have for some constant $C_l$ and for all $i$,
\[
\inf_{|z| \le K_X + K_0 } \ddot l_i (z) \ge2/ C_l .
\]
\end{condition}

\begin{example}
 The least absolute deviations loss is
$\rho(y,z) := | y-z| $. Let $G_i$ be distribution function of $Y_i$ given
$x_i$ ($i=1 , \ldots, n $). Then $f_i^0$ is the median of $G_i$, and
Condition~\ref{conB} requires that
$G_i$ has a strictly positive density $g_i$ on $\{ |z| \le K_X + K_0 \}$
for all $i$.
\end{example}

We now define
\[
\Gamma(S_0) := 16 C_l  \biggl[ { s_0 \over\phi^2 (3, S_0)}  \biggr]
.
\]
Fix some $t>0$ and define
\[
\lambda_{\varepsilon} (t) :=
16 K_X \sqrt{2(t + \log p) \over n}.
\]
The following theorem is a reformulation of results~in \citet
{vandeG07}, \citet{vandeG07} or \citet{BvdG2011}.

\begin{theorem} \label{oracle-robust.theorem}
Let $\hat\beta$ be the $\ell_1$-penalized robust estimator.
Assume Conditions~\ref{conA1},~\ref{conA2} and
\ref{conB}. Suppose that
%
\begin{equation}\label{s02}
\lambda_{\varepsilon} (t) \Gamma_0 (S_0) \le\tfrac{1}{4} .
\end{equation}
Take
\[
4 \lambda_{\varepsilon} (t) \le\lambda\le\frac{ 1}{\Gamma(S_0) } .
\]
With probability at least $1- \alpha(t)$, where
$\alpha(t) :=\break 3 \exp[-t] $,
it holds that
\[
\| \hat\beta- \beta^0 \|_1 \le\frac{\lambda}{2} \Gamma(S_0)
\]
and
\[
\| \hat f_{\hat\beta} - f_{\beta^0} \|_n^2 \le
\tfrac{3}{4} C_l \lambda^2 \Gamma(S_0) .
\]
\end{theorem}

\begin{remark}\label{oracle-robust.remark}
Similar remarks can be made as for the $\ell_1$-penalized
quasi-likelihood estimator. The new element in the result is that with
robustness the tuning parameter
$\lambda$ does not depend on some noise level $\sigma$.
\end{remark}

\section{Variable Selection with Quasi-Likelihood Loss}
\label{selection-quasi-likelihood.section}

Note that the bounds for the $\ell_1$-error $\| \hat\beta- \beta^0\|_1$,
given in Theorems~\ref{oracle-quasi-likelihood.theorem}
and~\ref{oracle-robust.theorem}, can be invoked to show that,
with large probability, the $\ell_1$-regularized estimator
will detect most of the nonzero coefficients $\beta^0 $ which are
large enough:
for all $\eta>0$,
\begin{eqnarray*}
&&\#\{ \hat\beta_j \not= 0 ,  | \beta_j^0|\ge\lambda/ \eta \} \\
&&\quad
\ge
\# \{ |\beta_j^0 | \ge\lambda/\eta\} - \eta\| \hat\beta- \beta
^0 \|_1 / \lambda.
\end{eqnarray*}
In other words, if a large proportion of the nonzero coefficients is
sufficiently far above the noise
level in absolute value, then
there will also be many true positives.
By this argument, if all nonzero coefficients of $\beta^0$ are of
larger order than $\lambda\Gamma(S_0) $,
we will have
$\hat S \supset S_0$, where
\[
\hat S := \{ j \dvtx  \hat\beta_j \not= 0 \} .
\]

This section will study the false positives. We show that for the case
of quasi-likelihood loss,
an irrepresentable condition similar to Meinshausen and \linebreak[4] B{\"u}hlmann (\citeyear{mebu06}) and \citet{ZY07}
implies that there
are no false positives, that is, that $\hat S \subset S_0$. Such result
can also be obtained for robust
loss, but is omitted here.

\subsection{The Case of Least Squares Loss}
Again, as preparation, let us first consider the standard linear model
and the least squares Lasso estimator,
\[
\hat\beta= \operatorname{arg}\mathop{\min}_{\beta\in\R^p }  \{ \| Y- X \beta\|
_n^2 + \lambda\| \beta\|_1  \} .
\]
Let $\mathbf{X} (S) := ( \mathbf{X}_j )_{j \in S}$ be the design matrix
consisting of the variables
in $S$, and let
\begin{eqnarray*}
\hat\Sigma_{1,1} (S) &:=& \mathbf{X}^T (S) \mathbf{X} (S) / n ,\\
  \hat\Sigma
_{1,2} (S) &:=& \mathbf{X}^T (S^c) \mathbf{X} (S) / n .
\end{eqnarray*}
In \citet{BvdG2011} (Exercise 7.5) or \citet{Geer:11}, one can find the
following result.

\begin{theorem}\label{teo7.1}
Suppose that $\lambda>\lambda_0 $ where
$\lambda_0 \ge2 \max_{1 \le j \le p } |\varepsilon^T \mathbf{X}_j |/ n $.
Assume moreover the
irrepresentable condition
\[
\sup_{\| \tau_{S_0} \|_{\infty} \le1 } \| \hat\Sigma_{2,1} (S_0)
\hat\Sigma_{1,1}^{-1} (S_0)
\tau_{S_0} \|_{\infty} < \frac{ \lambda- \lambda_0}{\lambda+
\lambda_0 } .
\]
Then $\hat S \subset S_0$.
\end{theorem}

We remark that an irrepresentable condition (see also below in
Definition~\ref{irrepresentable.definition}) is always rather strong.
However, for
exact variable selection, an irrepresent\-able condition
is essentially necessary,
as shown in \citet{mebu06}, \citet{ZY07}, \citet{BvdG2011}. By
thresholding the estimated coefficients
and refitting, or by applying the adaptive Lasso, one can
often improve on
variable selection and yet maintain a good prediction and estimation error.
The conditions for the latter are much less restrictive than the irrepresentable
condition. We refer to \citet{Geer:11} for details.

\subsection{General Quasi-Likelihood Loss}
The results are based on
he Karush--Kuhn--Tucker (or KKT-) conditions; see
\citet{bertsimas1997introduction}.
In our context, they read as follows:

\begin{KKTConditions*}
We have
\[
\frac{ \partial}{\partial\beta}
\frac{1}{n } \sum_{i=1}^n Q (Y_i , x_i^T \beta)  \Bigg\vert_{\beta
= \hat\beta}
=-\lambda\hat\tau.
\]
Here $\|\hat \tau\|_{\infty} \le1$, and moreover
\[
\hat\tau_{j}\mathrm{l} \{ \hat\beta_{ j} \not= 0 \} = \mathrm{sign}(
\hat\beta_{j} ) ,\quad  j = 1 , \ldots, p.
\]
\end{KKTConditions*}

Let
\[
\hat\Sigma_{j,k} := {1 \over n} \sum_{i=1}^n x_{i,j} x_{i,k}
w_{i}^2 ,
\]
where
\[
w_i^2 := h^2 ( x_i^T \beta^0) V\circ G( x_i^T \beta^0) ,\quad  i=1 ,
\ldots, n .
\]
Thus, $\hat\Sigma$ is the weighted Gram matrix
\[
\hat\Sigma= \mathbf{X}^T W^2 \mathbf{X} /n ,\quad  W^2 := \mathrm{diag} (w_1^2 ,
\ldots, w_n^2 ) .
\]
We write $\mathbf{X}_W := W \mathbf{X} $,
so that $\hat\Sigma= \mathbf{X}_W^T \mathbf{X}_W / n $.

Let $\mathbf{X}_W (S) $ be the weighted design matrix consisting of the
variables in $S$, and
\begin{eqnarray*}
\hat\Sigma_{1,1} (S)&:=& \mathbf{X}_W^T (S) \mathbf{X}_W (S)/n ,\\
\hat\Sigma_{2,1} (S) &:=& \mathbf{X}_W^T (S^c) \mathbf{X}_W (S)/n .
\end{eqnarray*}

\begin{definition} \label{irrepresentable.definition}
Let $0 < \theta\le1$ be given. We say that the $\theta$-irrepresentable condition
is met for the set $S$ if
\[
\max_{\| \tau_S \|_{\infty} \le1}
\|\hat \Sigma_{2,1} (S) \hat\Sigma_{1,1}^{-1} (S) \tau_S \|
_{\infty} \le\theta
.
\]
\end{definition}

Here is how the $\theta$-irrepresentable condition can be linked
with variable selection.

\begin{theorem} \label{select.theorem}
Let $0 \le\lambda_0 < \lambda$.
Suppose that
%
\begin{equation}\label{KKT2}
\hat\Sigma( \hat\beta- \beta^0 ) = - v ,
\end{equation}
where $ | v_j | \le\lambda+ \lambda_0$, and $v_j
\hat\beta_j \ge(\lambda- \lambda_0) |\hat\beta_j| $, $j=1 ,
\ldots, p$.
Suppose moreover the $\theta$-irrepresentable condition is met for $S_0$,
with $\theta< (\lambda- \lambda_0) / (\lambda+ \lambda_0 ) $.
Then $\hat S \subset S_0$.
\end{theorem}

In the proof of Theorem~\ref{selection-quasi-likelihood.theorem}
below, we show that
equation
(\ref{KKT2}) in
Theorem~\ref{select.theorem} holds for some $v$ satisfying the
conditions of
this theorem. This allows us then to conclude that $\hat S \subset S_0$.

As one sees in the KKT conditions, the derivative at $\hat\beta$ of
the loss function occurs.
We will need to compare this by the derivative at $\beta^0$. To bring
this to an end
we need,
in addition to Conditions~\ref{conA3} and~\ref{conA4}, certain Lipschitz conditions
on $h$ and $g$.

\renewcommand{\thecondition}{\textup{A}\arabic{condition}}
\setcounter{condition}{4}
\begin{condition} \label{conA5}
For $K_X$ and $K_0$, given in Conditions~\ref{conA1} and
\ref{conA2} respectively, we have for all $| z_0 | \le|z| \le K_X + K_0 $, and
some constant $L_h$,
\[
| h(z) - h( z_0) | \le L_h | z - z_0 | .
\]
\end{condition}

\begin{condition} \label{conA6}
For $K_X$ and $K_0$ given in Conditions~\ref{conA1} and
\ref{conA2} respectively, we have for all $| z_0 | \le|z| \le K_X + K_0 $, and
some constant $L_g$,
\[
|g(z) - g(z_0) | \le L_g | z - z_0 |/2 .
\]
\end{condition}

\begin{remark}$\!\!\!\!\!\!$
Under the additional Conditions~\ref{conA5} and~\ref{conA6}, one can
improve the
constants in Theorem~\ref{oracle-quasi-likelihood.theorem}. It is also clear
that Conditions~\ref{conA5} and~\ref{conA6} hold for least squares and logistic loss.
\end{remark}

With these new constants, we define
\[
L_{h,V}:=( L_g + L_h C_V) C_h ,\quad  L_{h,X}+ 16 L_h K_X^2 .
\]
We moreover let
\[
\Gamma_{\varepsilon} := \Gamma(S_0) := 16 C_{h,V} \Gamma_{\mathrm{effective}} (S_0) ,
\]
and
\[
\Gamma_0 := \Gamma_0 (S_0) := 6 L_{h,V} C_{h,V}^2 \Gamma_{\mathrm{effective}} (S_0) .
\]
Fix some $t>0$, and define
\[
\lambda_{\varepsilon} (t) :=
C_{h,X} \sigma\sqrt{\frac{2(t + \log p)}{n}}
\]
and
\[
\lambda_0 (t):=
L_{h,X} \sigma\sqrt{\frac{2( t + 2 \log p)}{n}}
.
\]
Define
\[
\alpha(t) := 9\exp[-t] + 9 \kappa^4 /(n \sigma^4).
\]
Thus, up to constants, $\Gamma_{\varepsilon}$ and $\Gamma_0$ are the effective
sparsity. Moreover, for $t \asymp\log n $ (say), $\lambda^{\varepsilon}
(t) \asymp\lambda_0 (t)
\asymp\sqrt{\log(p \vee n ) / n } $ and $\alpha(t) \asymp1/n $.

We arrive at the main result of this section.

\begin{theorem}\label{selection-quasi-likelihood.theorem}
Let $\hat\beta$ be the
$\ell_1$-penalized quasi-likelihood estimator.
Assume Conditions
\ref{conA} and
\ref{conA1}--\ref{conA6}.
Assume that (\ref{s0}) holds, that is,
\[
\lambda_{\varepsilon} (t) \Gamma_{\varepsilon} \le\gamma_1 \le \tfrac{1
}{4},
\]
where $\gamma_1$ is given by
\[
\gamma_1 := \frac{ \lambda_{\varepsilon} (t)}{\lambda} .
\]
Assume now that
%
\begin{equation} \label{s03}
\lambda_{\varepsilon} (t) \Gamma_0 \le\gamma_1 \gamma_{\varepsilon}\quad
\mbox{for  some }  \gamma_{\varepsilon} <
1- \gamma_1 ,
\end{equation}
as well as
%
\begin{equation}\label{s04}
\quad\lambda_0 (t) \Gamma_{\varepsilon} \le\gamma_0 \quad \mbox{for  some }
\gamma_0 < 1-
\gamma_{\varepsilon} - \gamma_1.
\end{equation}

Assume furthermore the $\theta$-irrepresentable condition with
\[
\theta< \frac{1- \gamma}{1+\gamma},\quad  \gamma:= \gamma_{\varepsilon} +
\gamma_0 + \gamma_1 .
\]

With probability at least $1- \alpha(t) $, it holds that
\mbox{$\hat S \subset S_0$}.
\end{theorem}

\begin{remark} \label{var.remark}
Let us take $\lambda_{\varepsilon} (t) \asymp\lambda_0 (t) \asymp
\lambda\asymp\sqrt{\log p /n } $. The constants $\gamma_0$,
$\gamma_1$ and $\gamma_{\varepsilon}$ are small,
depending on the constants appearing in Conditions~\ref{conA}
and~\ref{conA1}--\ref{conA6}. Fixing these, they can be kept away
from zero, and hence also the $\theta$-irrepresentable condition
is assumed for a value of $\theta$ that stays away from zero.
Conditions (\ref{s03}), and (\ref{s04}) again require that the effective
sparsity is
sufficiently smaller that $\sqrt{\log p / n }$. Formulated differently,
the results of Theorems~\ref{selection-quasi-likelihood.theorem} and~\ref{oracle-quasi-likelihood.theorem} imply that
if the $\theta$-irrepresentable condition holds, and if
$\Gamma_\mathrm{effective} (S_0) \le\gamma\sqrt{\log p / n }$ for
sufficiently small values of $\theta$
and $\gamma$
(depending only on the constants appearing in
Conditions
\ref{conA} and \mbox{\ref{conA1}--\ref{conA6}}), then with an appropriate choice of
$\lambda\asymp\sqrt{\log p / n } $ the Lasso estimator has with
large probability prediction error
$ \Gamma_\mathrm{effective} (S_0) \log p / n $, $\ell_1$-error
$\Gamma_\mathrm{effective} (S_0) \sqrt{\log p / n }$ and
no false positives.
\end{remark}

\section{Random Design}\label{random-design.section}
Consider quasi-likelihood loss.
It is easy to see that under the conditions of Theorem \ref
{oracle-quasi-likelihood.theorem},
one has with large probability
\[
( \hat\beta- \beta^0 )^T \hat\Sigma(\hat\beta- \beta^0) \le
6 C_{h,V}^3 \lambda^2 \Gamma_\mathrm{effective} (S_0).
\]
This follows from $w_{i}^2 \le C_{h,V} /2 $, where
as in Section~\ref{selection-quasi-likelihood.section},
$w_{i}^2 = h^2 (x_i^T \beta^0) V \circ G( x_i^T \beta^0)$,
$i=1 , \ldots, n$. Let $\Sigma$ be some other $p \times p$
positive semi-definite matrix. Then
\[
\| ( \hat\Sigma- \Sigma) (\hat\beta- \beta^0 )\|_{\infty} \le
\lambda_X \| \hat\beta- \beta^0 \|_1 ,
\]
where
\[
\lambda_X := \max_{j,k} | \hat\Sigma_{j,k} - \Sigma_{j,k} | .
\]
Thus, under the conditions of Theorem~\ref{oracle-quasi-likelihood.theorem},
one has that with large probability
\[
\| ( \hat\Sigma- \Sigma) (\hat\beta- \beta^0 )\|_{\infty} \le
\lambda\lambda_X \Gamma(S_0)/2 .
\]
One can verify that if $\lambda_X \Gamma(S_0)$ is small enough,
say for some $\gamma_X$ sufficiently small
%
\begin{equation}\label{s05}
\lambda_X \Gamma(S_0) \le\gamma_X ,
\end{equation}
then one may reformulate the compatibility condition
replacing $\| f_{\beta} \|_n^2$ by $\beta^T\Sigma\beta$, and the
theory for prediction
and $\ell_1$-error goes
through essentially without new arguments. One can then also establish bounds
for $(\hat\beta- \beta^0 )^T\Sigma(\hat\beta- \beta^0 )$.
Similarly, one may reformulate the $\theta$-irrepresentable condition
with $\hat\Sigma$ replaced by $\Sigma$, and obtain variable selection
without needing new arguments. In the case where $\Sigma$ is the
population version of $\hat\Sigma$, the latter built from an i.i.d. sample
of covariables, one can show that with large probability
$\lambda_X$ is of order $\sqrt{\log p / n }$. In other words (and modulo
the compatibility constant), then
condition (\ref{s05}) is another instance where it is required that
the sparsity $s_0$ is not of larger order than\break $\sqrt{n / \log p } $.
We refer to \citet{BvdG2011} for more precise statements.

\section{Conclusion}\label{discussion}
The results of this paper show that the oracle and variable
selection properties of the Lasso for the linear model
also hold for the generalized linear model. We prove this
under the assumption that the is sparsity sufficiently smaller than
$\sqrt{n / \log p }$.
We note that the results rely heavily on the convexity of the
loss function. This allows one to work with
an unbounded parameter space. If the estimators are a priori
restricted to lie in a given bounded set, one can extend
the results to nonconvex loss [see \citet{Staedler:10} for the
mixture model, and
\citet{schelldorfer2011estimation} for the mixed effects model]
and one can
moreover prove oracle results for the almost
linear in $s_0$ regime of sparsity.

\begin{supplement}[id=supp]
\stitle{Supplementary material for ``Quasi-likelihood and/\break or robust estimation in high dimensions''\\}
\slink[doi]{10.1214/12-STS397SUPP}
\sdatatype{.pdf}
\sfilename{STS397\_supp.pdf}
\sdescription{Due to space constraints, the proofs and technical
details have been given
in the supplementary document \protect\citet{GeerMullersupp:12}.}
\end{supplement}

\section*{Acknowledgments}
Research supported by SNF 20PA21-120050.

%

\end{document}